\documentclass{article}
\usepackage{amsmath, amssymb, amsthm, relsize, tikz-cd}

\theoremstyle{plain}
\newtheorem{theorem}{Theorem}[section]

\theoremstyle{definition}

\DeclareMathOperator*{\OmSum}{\mathlarger{\mathlarger{\Omega}}}

\newcommand{\T}{\mathcal{T}}

\begin{document}

\title{Second Order Transfer Equations; and Generalizations to Arbitrary Orders}

\author{James David Nixon\\
	JmsNxn92@gmail.com\\}

\maketitle

\begin{abstract}
The author provides a solution to the equation $ y(s+2) = \T^2 y = F(s,y,\T y) = F(s,y(s),y(s+1))$; where $y$ is holomorphic; and $F$ is a holomorphic function with specific decay conditions. This result is provided using infinite compositions, and a limiting process. The technique is generalized to arbitrary $k$'th order transfer equations: $u(s+k) = F(s,u(s),u(s+1),...,u(s+k-1))$. The technique is derived by utilizing solutions $w(s+1) = \T w = F(s,w)$ and sequentially approximating $u$, or $y$, with a sequence of said $w$.
\end{abstract}



\section{Introduction}\label{sec1}
\setcounter{equation}{0}

To paint the picture of what we're attempting in this paper is a nuanced thing to do. We can point expressly to our previous work in \cite{Nix} which aimed at solving First Order Transfer Equations and First Order Difference Equations in the complex plane. Of this, we change our focus to Second Order Transfer Equations. To give a synopsis of what this problem is, it is helpful to dive right in.

Suppose $F(s,z,w)$ is a holomorphic function taking $\mathcal{S} \times \mathcal{G} \times \mathcal{G} \to \mathcal{G}$, where $\mathcal{S}$ and $\mathcal{G}$ are domains in the complex plane. Is it possible to find a holomorphic function $y:\mathcal{S} \to \mathcal{G}$ such that,

\[
y(s+2) = F(s,y(s),y(s+1))?
\]

One can readily see the connection to second order differential equations $y'' = f(s,y,y')$. Where, now the operator $\frac{d}{ds}$ is replaced with $\mathcal{T}y = y(s+1)$. A similar idea which requires a very different tool-set. We are trying to solve,

\[
\T^2 y = F(s,y, \T y)\\
\]

And not only that; we want to solve this analytically in $\mathbb{C}$. Or at least, within domains of $\mathbb{C}$. We'll deal entirely with holomorphic functions; which, although appearing more difficult to solve, the discussion reserved is rather concise.

This type of problem can be traced back to Leonardo Fibonacci. If $F(s,z,w) = z + w$, then the Fibonacci sequence is a solution on $\mathbb{N}$, and its solutions in the complex plane are given by Binet,

\[
y(s) = \frac{\phi^s - \psi^{s}}{\phi - \psi}
\]

where $\phi = \frac{1 + \sqrt{5}}{2}$ and $\psi = \frac{1 - \sqrt{5}}{2}$, and one can choose from all branches of the exponential function. Which is to say,

\[
y(s+2) = y(s) + y(s+1) = F(s,y(s),y(s+1))\\
\]

When $s$ is restricted to the natural numbers, these recursions generate a sequence; so long as we are given some initial conditions. When $s$ is allowed to exist in the complex plane the problem gains difficulty. There is no natural way to obtain the value $s \in \mathbb{C}$. All that can be said is once a value at $s_0-1$ and $s_0-2$ is obtained; the value at $s_0$ can be found. Then the value at $s_0+1$ can be found, so on and so forth.

Which, to our next point; on the initial value problem of Second Order Recursions, there's very little we can say. We can say there are solutions dropping the initial value problem, and that's about it. Similarly, the uniqueness of these solutions is a mystery we daren't attempt to attack.\\

In order to broach this problem we will need the language of infinite compositions. Sadly there exists no crash course in the language; or the theory. If the reader wishes for a more thorough discussion of the notation and the language, we refer to \cite{Nix, Nix2}. The author will simply state the definition of our object of interest. If $h_j(s,z,w)$ is a sequence of analytic functions, the expression

\[
\OmSum_{j=1}^n h_j(s,z,w) \bullet w = h_1(s,z,h_2(s,z,...h_n(s,z,w)))
\]

The bullet implies which variable we are composing across, and the symbol $\OmSum$ behaves no differently than Euler's usage of $\sum$ and $\prod$.

The case will be taken that when either of the variable is fixed to $z=A$ or $w=A'$, which through linear substitution, can be reduced to the case $z=0$ or $w=0$. The author has developed a theorem in \cite{Nix2} which was stated without generality there, but will be stated in full generality here. The proof can be found in the aforementioned paper, though worded in a restricted case; we will addend it here in the next section.

\begin{theorem}\label{thmA}
Let $\{H_j(s,z)\}_{j=1}^\infty$ be a sequence of holomorphic functions such that $H_j(s,z) : \mathcal{S} \times \mathcal{G} \to \mathcal{G}$ where $\mathcal{S}$ and $\mathcal{G}$ are domains in $\mathbb{C}$. Suppose there exists some $A \in \mathcal{G}$, such for all compact sets $\mathcal{N}\subset\mathcal{G}$, the following sum converges,

\[
\sum_{j=1}^\infty ||H_j(s,z) - A||_{z \in \mathcal{N},s \in \mathcal{S}} < \infty
\]

Then the expression,

\[
H(s) = \lim_{n\to\infty}\OmSum_{j=1}^n H_j(s,z)\bullet z = \lim_{n\to\infty} H_1(s,H_2(s,...H_n(s,z)))\\
\]

Converges uniformly for $s \in \mathcal{S}$ and $z \in \mathcal{N}$ as $n\to\infty$ to $H$, a holomorphic function in $s\in\mathcal{S}$, constant in $z$.
\end{theorem}

This theorem alone allows for the solution of $k$'th Order Transfer Equations, but we will lead by example with Second Order Transfer Equations. In this case we can find the philosophy for the $k$'th case, albeit subtracting a whole mess of symbols which arise for arbitrary orders. Imagine this paper as a proof structure using the example of $2\times2$ matrices; while maintaining the abstract language necessary to resolve $k\times k$ matrices.

\section{Proof Of Theorem \ref{thmA}}

\begin{proof}

The first thing we show is for all $\epsilon > 0$, there exists some $N$, such when $m \ge n  > N$,

\[
|\OmSum_{j=n}^{m} H_j(s,z)\bullet z - A| < \epsilon
\]

For $z$ in $\mathcal{N}\subset \mathcal{G}$ (where $A$ is in the open component of $\mathcal{N}$), and $s\in\mathcal{S}$. This then implies as we let $m\to\infty$, the tail of the infinite composition stays bounded. Forthwith, the infinite composition becomes a normal family, and proving convergence becomes simpler. We provide a quick proof of this inequality.\\

Set $||H_j (s, z)-A||_{\mathcal{S}, \mathcal{N}} = \rho_j$. Pick $\epsilon > 0$, and choose $N$ large enough so when $n > N$,

\[
\rho_n < \epsilon
\]

Denote: $\phi_{nm}(s, z) =\OmSum_{j=n}^m H_j (s, z) \bullet z = H_n(s, H_{n+1}(s, ...H_m(s, z)))$. We go by induction on the difference $m-n = k$. When $k=0$ then,

\[
||\phi_{nn}(s,z) - A||_{\mathcal{S},\mathcal{N}} = ||H_n(s,z)-A||_{\mathcal{S},\mathcal{N}}= \rho_n < \epsilon
\]

Assume the result holds for $m-n < k$, we show it holds for $m-n = k$. Observe,

\begin{eqnarray*}
||\phi_{nm}(s,z)-A||_{\mathcal{S},\mathcal{N}} &=& ||H_n(s,\phi_{(n+1)m}(s,z)) - A||_{\mathcal{S},\mathcal{N}}\\
&\le& ||H_n(s,z)-A||_{\mathcal{S},\mathcal{N}}\\
&=& \rho_n < \epsilon\\
\end{eqnarray*}

Which follows by the induction hypothesis because $\phi_{(n+1)m}(s,z) \subset \mathcal{N}$--it's in a neighborhood of $A$ which is in $\mathcal{N}$. That is $m-n-1 < k$.

The next step is to observe that $\OmSum_{j=1}^m H_j(s,z)$ is a normal family as $m\to\infty$, for $z \in \mathcal{N}$ and $s \in \mathcal{S}$. This follows because the tail of this composition is bounded. Therefore we can say $||\OmSum_{j=1}^m H_j(s,z)||_{\mathcal{S},\mathcal{N}} < M$ for all $m$.

Since $\phi_m(s,z) = \OmSum_{j=1}^m H_j(s,z)\bullet z$ are a normal family for all compact sets $\mathcal{N}\subset \mathcal{G}$; there is some constant $M \in \mathbb{R}^+$ and $L \in \mathbb{R}^+$ such,

\[
||\frac{d^k}{dz^k} \phi_m(s,z) ||_{\mathcal{S},\mathcal{N}} \le M \cdot k! \cdot L^k
\]

(Use Cauchy's Integral Theorem.) Secondly, using Taylor's theorem,

\begin{eqnarray*}
\phi_{m+1}(s,z) - \phi_m(s,z) &=& \phi_m(s,H_{m+1}(s,z)) - \phi_m(s,z)\\
&=& \sum_{k=1}^\infty \frac{d^k}{dz^k} \phi_m(s,z) \frac{(H_{m+1}(s,z) - z)^k}{k!}\\
&=& (H_{m+1}(s,z) - z) \sum_{k=1}^\infty \frac{d^k}{dz^k} \phi_m(s,z) \frac{(H_{m+1}(s,z) - z)^{k-1}}{k!}\\
\end{eqnarray*}

So that, setting $z=A$,

\begin{eqnarray*}
||\phi_{m+1}(s,A) - \phi_m(s,A)||_{s \in\mathcal{S}} &\le& ||H_{m+1}(s,A) - A||_{s\in\mathcal{S}} \sum_{k=1}^\infty M L^k ||H_{m+1}(s,A) - A||^{k-1}\\
&\le& ||H_{m+1}(s,A) - A||_{\mathcal{S}} \frac{ML}{1-q}\\
\end{eqnarray*}

For $L||H_{m+1}(s,A) - A||_{\mathcal{S}} \le q <1 $, which is true for large enough $m>N$. Setting $C = \frac{ML}{1-q}$. Applying from here,

\[
||\phi_{m+1}(s,A) - \phi_m(s,A)||_{s \in \mathcal{S}} \le C ||H_{m+1}(s,A) - A||_{s \in \mathcal{S}}\\
\]

This is a convergent series per our assumption. Choose $N$ large enough, so that when $m,n>N$,

\[
\sum_{j=n}^{m-1}||H_{j+1}(s,A) - A||_{s \in \mathcal{S}} < \frac{\epsilon}{C}\\
\]

Then,

\begin{eqnarray*}
||\phi_{m}(s,A) - \phi_n(s,A)||_{s \in \mathcal{S}} &\le& \sum_{j=n}^{m-1} ||\phi_{j+1}(s,A) - \phi_j(s,A)||_{s \in \mathcal{S}}\\
&\le& C\sum_{j=n}^{m-1}||H_{j+1}(s,A) - A||_{s \in \mathcal{S}}\\
&<& \epsilon
\end{eqnarray*}

So we can see $\phi_m(s)$ must be uniformly convergent for $s \in \mathcal{S}$, and therefore defines a holomorphic function $H(s)$ as $m\to\infty$.

This tells us,

\[
H(s) = \OmSum_{j=1}^\infty H_j(s,z)\bullet z \Big{|}_{z=A}\\
\]

Converges and is holomorphic. To show this function equals,

\[
\OmSum_{j=1}^\infty H_j(s,z)\bullet z
\]

For all $z \in \mathcal{G}$; simply notice that,

\[
\OmSum_{j=m}^\infty H_j(s,z)\bullet z
\]

Is arbitrarily close to $A$ as we let $m$ grow (which was shown at the beginning of this proof). Then,

\begin{eqnarray*}
\OmSum_{j=1}^\infty H_j(s,z)\bullet z &=& \OmSum_{j=1}^{m-1} H_j(s,z)\bullet \OmSum_{j=m}^\infty H_j(s,z)\bullet z\\
&=& \lim_{m\to\infty} \OmSum_{j=1}^{m-1} H_j(s,z)\bullet \lim_{m\to\infty} \OmSum_{j=m}^\infty H_j(s,z)\bullet z\\
&=& \OmSum_{j=1}^\infty H_j(s,z)\bullet z\Big{|}_{z=A}\\
\end{eqnarray*}

This concludes our proof.

\end{proof}

\section{Constructing a sequence}\label{sec2}
\setcounter{equation}{0}

Let $\mathcal{S}$ be a strip in the complex plane, $\mathcal{S}= \{s\in\mathbb{C}\,:\,|\Im(s)| < a\}$. Let $\mathcal{G}$ be a domain in $\mathbb{C}$. Let $F(s,z,w)$ be a holomorphic function for $s \in \mathcal{S}$, and $z,w \in \mathcal{G}$. Further, $F$ lives in $\mathcal{G}$; $F(s,z,w)\in \mathcal{G}$. Let's also assume the following sum converges,

\[
\sum_{j=1}^\infty ||F(s-j,z,w)||_{\mathcal{S}_J, \mathcal{N},\mathcal{N}} < \infty
\]

for all compact sets $\mathcal{N}\subset \mathcal{G}$; and $\mathcal{S}_J = \{s \in \mathcal{S}\,:\,\Re(s) < -J,\, J \in \mathbb{R}^+\}$, and $||...||$ denotes the supremum norm. By Theorem \ref{thmA}, the function,

\[
y_0(s,w) = \OmSum_{j=1}^\infty F(s-2j,z,w) \bullet z
\]

Is holomorphic for $s \in \mathcal{S}_J$ and $w \in \mathcal{G}$. It also satisfies the functional equation,

\[
y_0(s+2,w) = F(s,y_0(s,w),w)
\]

This function equally has decay to $0$ as $s \to - \infty$. From this identification we can say that,

\[
\sum_{j=1}^\infty ||F(s-2j,z,y_0(s+1-2j))||_{\mathcal{S}_J, \mathcal{N}} < \infty
\]

And by inducting this process, we can define a sequence,

\[
y_n(s) = \OmSum_{j=1}^\infty F(s-2j,z,y_{n-1}(s+1-2j)) \bullet z
\]

Such that $y_n(s)$ is holomorphic for $s \in \mathcal{S}_J$, and satisfies the functional equation,

\[
y_n(s+2) = F(s,y_n(s),y_{n-1}(s+1))
\]

If this were to converge as $n\to\infty$, we would have an existence and construction method for our desired function $y$. In the next section we provide a proof of said fact.

\section{Showing convergence of $y_n$}\label{sec3}
\setcounter{equation}{0}

Take $J$ very large. Assume $\Re (s) < -J$. It is no hard task to show we can choose $J$ large enough such that there is some $\lambda<1/2 $, where for all $s \in \mathcal{S}_J = \{s \in \mathcal{S}\,:\,\Re (s) < -J\}$,

\[
||F(s,z,w) - F(s,z,w')||_{\mathcal{S}_J,z \in \mathcal{N}} \le \lambda | w - w'|
\]

And similarly for $z$,

\[
||F(s-j,z,w) - F(s-j,z',w)||_{\mathcal{S}_J,w \in \mathcal{N}} \le \lambda | z - z'| 
\]

And just as well, we have,

\[
||F(s-j,z,w) - F(s-j,z',w')||_{\mathcal{S}_J} \le \lambda | z - z'| + \lambda | w - w'| 
\]

Looking again at the functional equation $y_n$ exhibits;

\[
y_n(s) = F(s-2,y_n(s-2),y_{n-1}(s-1))
\]

We get,

\begin{eqnarray*}
|y_n(s) - y_{n-1}(s)| &\le& \lambda | y_n(s-2) - y_{n-1}(s-2)|\\ &+& \lambda | y_{n-1}(s-1) - y_{n-2}(s-1)| \\
\end{eqnarray*}

This inequality alone justifies convergence of the sequence. If we iterate the procedure; recursively expanding $N$ times,

\begin{eqnarray*}
|y_n(s) - y_{n-1}(s)| &\le& \lambda^N | y_n(s-2N) - y_{n-1}(s-2N)|\\ &+& \sum_{k=1}^{N} \lambda^k | y_{n-1}(s-2k+1) - y_{n-2}(s-2k+1)|\\
\end{eqnarray*}

The first term can be made arbitrarily small for large enough $N$ while $s \in \mathcal{S}_J$. Let $N \to \infty$. Taking the supremum norm across $\mathcal{S}_J$, and bounding each $|y_{n-1}(s-2k+1) - y_{n-2}(s-2k+1)| \le || y_{n-1}(s) - y_{n-2}(s)||_{\mathcal{S}_J}$; the second term can be bounded as,

\[
||y_n(s) - y_{n-1}(s)||_{\mathcal{S}_J} \le \frac{\lambda}{1-\lambda} || y_{n-1}(s) - y_{n-2}(s)||_{\mathcal{S}_J}\\
\]

Now $\lambda < 1/2$, and therefore $\frac{\lambda}{1-\lambda} = \mu < 1$. Let $A = ||y_1(s) - y_0||_{\mathcal{S}_J}$, and take $n,m >N$ so that $\sum_{k=n}^{m-1} \mu^k < \dfrac{\epsilon}{A}$. By using a telescoping series;

\begin{eqnarray*}
||y_n(s)-y_m(s)||_{\mathcal{S}_J} &\le& \sum_{k=n}^{m-1} ||y_{k+1}(s) - y_k(s)||_{\mathcal{S}_J}\\
&\le& \sum_{k=n}^{m-1} \mu^k ||y_1 - y_0||_{\mathcal{S}_J}\\
&\le& \sum_{k=n}^{m-1} A\mu^k\\
&<& \epsilon\\
\end{eqnarray*}

Therefore the limit converges, and $y_n(s) \to y(s)$ uniformly while $s \in \mathcal{S}_J$. We can extend $y$, due to the functional equation, to all of $\mathcal{S}$:

\[
y(s+2) = F(s,y(s),y(s+1))\\
\]

So the domain of the function $y(s)$ is extendable to the entire strip $\mathcal{S}$. In order to extend this result to arbitrary $k$'th Order Recursions we simply need construct a function,

\[
y_n(s+k) = F(s,y_n(s),y_{n-1}(s+1), y_{n-1}(s+2),...,y_{n-1}(s+k-1))
\]

Which can be written as,

\[
y_n(s) = \OmSum_{j=1}^\infty F(s-jk,z_1,y_{n-1}(s-jk+1), y_{n-1}(s-jk+2),...,y_{n-1}(s-jk +k-1))\bullet z_1\\
\]

Carried from the initial condition,

\[
y_0(s) = \OmSum_{j=1}^\infty F(s-jk,z_1,z_2,...,z_k)\bullet z_1\\
\]

Similarly as before, when we take the difference, if $J$ is large enough, and $\lambda < 1/k$,

\[
||y_n - y_{n-1}|| \le \lambda||y_n-y_{n-1}|| + (k-1)\lambda||y_{n-1} - y_{n-2}||\\
\]

Upon,

\[
||y_n - y_{n-1}|| \le \lambda^N||y_n-y_{n-1}|| + (k-1) \sum_{j=1}^N \lambda^j||y_{n-1} - y_{n-2}||\\
\]

So that, for some $\mu < 1$ where $\mu = \frac{(k-1)\lambda}{1-\lambda}$,

\[
||y_n - y_{n-1}|| \le \mu ||y_{n-1} - y_{n-2}||\\
\]
And just as written above, the limit converges on $\mathcal{S}_J$ and we can extend to $\mathcal{S}$ in no more difficult a fashion. All of our discussions throughout this paper compress into the nice theorem:

\begin{theorem}[Existence Of Solutions For $k$'th Order Transfer Equations]
Let $\mathcal{S}$ be a horizontal strip in $\mathbb{C}$, and $\mathcal{G}\subseteq\mathbb{C}$ be a domain. Let $F(s,z_1,z_2,...,z_k): \mathcal{S} \times \mathcal{G}^k \to \mathcal{G}$ be a holomorphic function. Suppose for all compact sets $\mathcal{K}\subset\mathcal{S}$ and $\mathcal{N}\subset \mathcal{G}$,

\[
\sum_{j=1}^\infty ||F(s-j,z_1,z_2,...,z_k) - A||_{\mathcal{K},\mathcal{N}^k} < \infty\\
\]

For some $A \in \mathbb{C}$. Then there exists a holomorphic function $y:\mathcal{S}\to\mathcal{G}$ such that,

\[
y(s+k) = F(s,y(s),y(s+1),...,y(s+k-1))\\
\]
\end{theorem}

With this we move on to some examples before closing this brief notice.

\section*{Some Examples}

Here we provide some examples as to what kinds of $k$'th Order Transfer Equations we have provided solutions for. If we are given the function:

\[
F(s,z_1,z_2,...,z_k) = \sum_{j=1}^k e^{s+z_j}
\]

Then, we have constructed a solution $y(s) : \mathbb{C} \to \mathbb{C}$ such that,

\[
y(s+k) = \sum_{j=0}^{k-1} e^{s+y(s+j)}
\]

If we are given the function,

\[
F(s,z_1,z_2,...,z_k) = \frac{z_1 + z_2^2 + ... +z_k^k}{1+s^2}\\
\]

Then we have constructed a solution for $y(s) : \mathbb{C}/\{m\pm i\,:\, m \in \mathbb{N},\, m\ge 1\}\to\mathbb{C}$ such that,

\[
y(s+k) = \frac{y(s) + y(s+1)^2 + ... + y(s+k-1)^k}{1+s^2}\\
\]

Or maybe we are given,

\[
F(s,z_1,z_2, ...,z_k) = e^{z_1z_2\cdots z_k - s^2}\\
\]

Then we've constructed a solution $y:\mathbb{C} \to \mathbb{C}$,

\[
y(s+k) = e^{y(s)y(s+1)\cdots y(s+k-1) - s^2}\\
\]

Which is as nonsensical as it looks. The author can't even imagine what it looks like graphed. But $y$'s asymptotic to $e^{(s-k)^2}$ as $\Re(s) \to \pm \infty$, which is nice. It's in the Schwartz space. Though I doubt you can say much about the Fourier transform of it.

\section*{In Conclusion}

We thank the reader here for taking the time to read this brief. There exists many more problems on the horizon--relating to difference equations and techniques towards constructing solutions. And dare I say, it looks fruitful.

One such problem being, the initial condition problem. Now, if the reader has cared to notice, we have avoided all conversation of the extraneous variables in our construction. Which is to say, $y$ is not only holomorphic for $s$, by way of our construction, it is holomorphic for $z_2,z_3,..,z_k$ as well. Which, act as parameters, which in some way shape or form is responsible for our initial condition. That is, the value $y(0), y(1),y(2),...,y(k-1)$. A discussion of this, is definitely in order; but it's a problem the author hasn't concerned himself with.

Another such problem being, a clean way of taking transfer equations to difference equations. Of it, for instance, translating everything in this paper to the language of difference equations. And instead solving,

\[
y(s+2) - 2y(s+1) + y(s) = \Delta^2 y = G(s,y,\Delta y) = G(s,y(s),y(s+1)-y(s))\\
\]

And similarly for higher-order equations. Which the author assumes is some change of variables to go from $G$ to $F$ and back; but he lacks the mathematical gene which allows for such contemplation and result. Then the natural question is to talk about $\Delta_h y = \frac{y(s+h)-y(s)}{h}$, and limiting $h\to 0$; similarly for higher order equations. Then coming full circle to equations of the form,

\[
y^{(n)} = H(s,y,y',..,y^{(n-1)})\\
\]

Simply put, the author lacks the vocabulary for such an investigation. \textit{To mire, and to spit a spindle; so perhaps someone may correct--that is the prerogative.}

\end{document}